\documentclass[12pt,leqno]{article}
\usepackage{amsmath,amsthm}

\def\init{\setcounter{equation}{0}}
\newtheorem{theorem}{Theorem}[section]

\newcommand{\R}{{\bf R}}
\newcommand{\C}{{\bf C}}

\newtheorem{lemma}{Lemma}[section]

\newcommand{\e}{{\varepsilon}}

\title{Inverse problems for hyperbolic equations.
\author{G.Eskin, \ \ \  Department of Mathematics, UCLA,\\ Los Angeles,
CA 90095-1555, USA. \ E-mail: eskin@math.ucla.edu}
}

\begin{document}
\maketitle

\section{Formulation of the problem and the main theorem.}
\label{section 1}
\init

Let $\Omega$ be a smooth bounded domain in $\R^n, n\geq 2$.
Consider in the cylinder $\Omega\times (0,T_0)$ the following hyperbolic equation:

\begin{eqnarray}                               \label{eq:1.1}
Lu\stackrel{def}{=}
\left(-i\frac{\partial u}{\partial t}+A_0(x,t)\right)^2u(x,t)
\ \ \ \ \ \ \ \ \ \ \ \ \ \ \ \ \ \ \ \ \ \ \ \ \ \ \ \ \ \ 
\\
-\sum_{j,k=1}^n\frac{1}{\sqrt{g(x)}}
\left(-i\frac{\partial}{\partial x_j}+A_j(x,t)
\right)
 \sqrt{g(x)}g^{jk}(x)
\left(-i\frac{\partial}{\partial x_k}+A_k(x,t)\right)u
\nonumber
\\
-V(x,t)u=0,
\nonumber
\end{eqnarray}
where $A_j(x,t),0\leq j\leq n, V(x,t)$ are 
$C^\infty(\overline{\Omega}\times[0,T_0])$
functions,  analytic in $t, \|g^{jk}(x)\|^{-1}$ is the metric tensor
in $\overline{\Omega}, g(x)=\det \|g^{jk}\|^{-1}$.
We consider the initial-boundary value problem for (\ref{eq:1.1}) in
$\Omega\times(0,T_0)$:
\begin{equation}                                 \label{eq:1.2}
u(x,0)=u_t(x,0)=0, \ \ \ \ x\in \Omega,
\end{equation}
\begin{equation}                                 \label{eq:1.3}  
u(x,t)\left|_{\partial\Omega\times(0,T_0)}\right. = f(x,t).
\end{equation}

The following operator is called the Dirichlet-to-Neumann (D-to-N) operator:
\begin{equation}                                   \label{eq:1.4}
\Lambda f\stackrel{def}{=}\sum_{j,k=1}^n g^{jk}(x)\left(\frac{\partial u}{\partial x_j}
+iA_j(x,t)u\right)\nu_k\left
(\sum_{p,r=1}^ng^{pr}(x)\nu_p\nu_r\right)^{-\frac{1}{2}}
\left|_{\partial\Omega\times(0,T_0)}\right.,
\end{equation}
where $u(x,t)$ is the solution of the initial-boundary value problem
(\ref{eq:1.1}), (\ref{eq:1.2}), (\ref{eq:1.3}),  $\nu=(\nu_1,...,\nu_n)$
is the unit exterior normal vector at $x\in \partial\Omega$ with respect
to the Euclidian metric.  If $F(x)=0$ is the equation of $\partial\Omega$
in some neighborhood of a point  $x_0\in\partial\Omega$ then 
$\Lambda f$ has the following form in this neighborhood:
\begin{eqnarray}                                   \label{eq:1.5}
\Lambda f=\sum_{j,k=1}^n g^{jk}(x)\left(\frac{\partial u}{\partial x_j}
+iA_j(x,t)u\right)F_{x_j}(x)
\\
\cdot\left(\sum_{p,r=1}^ng^{pr}(x)F_{x_p}F_{x_r}
\right)^{-\frac{1}{2}}
\left|_{F(x)=0,0< t< T_0}\right. .
\nonumber
\end{eqnarray}
Let $\Gamma_0$ be an open subset of $\partial\Omega$.  We shall consider smooth
 $f(x,t)$ such that $\mbox{supp\ } f\subset\Gamma_0\times(0,T_0].$
The inverse problem consists of recovering the coefficients of (\ref{eq:1.1})
knowing the restriction of $\Lambda f$ to $\Gamma_0\times(0,T_0)$ for all
smooth $f$ with supports in $\Gamma_0\times(0,T_0]$.

There is a built-in nonuniqueness of this inverse problem:

a)  Let $y=\varphi(x)$ be a diffeomorphism of $\overline{\Omega}$ onto
$\overline{\Omega_0}\stackrel{def}{=}\varphi(\overline{\Omega})$ such
that $\Gamma_0\subset \partial\Omega_0$ and $\varphi=I$ on $\Gamma_0$.

Let $\hat{L}\hat{u}=0$ be the equation (\ref{eq:1.1}) in $y$-coordinates and
let $\hat{\Lambda}$ be the new D-to-N operator.  It follows from (\ref{eq:1.5})
that $\hat{\Lambda}=\Lambda$ on $\Gamma_0\times(0,T)$, i.e. 
$\hat{\Lambda}f|_{\Gamma_0\times(0,T_0)}=\Lambda f|_{\Gamma_0\times(0,T_0)}$
for all $f,\ \mbox{supp\ }f\subset\Gamma_0\times(0,T_0]$,  
i.e. the D-to-N operator on $\Gamma_0\times(0,T_0)$ cannot distinguish
between $Lu=0$ in $\Omega\times(0,T_0)$ and
$\hat{L}\hat{u}=0$ in $\Omega_0\times(0,T_0)$.

b)  Let $G_0(\Omega\times[0,T_0])$ be a group of 
$C^\infty(\overline{\Omega}\times[0,T_0])$ complex-valued functions 
$c(x,t)$ such that $c(x,t)\neq 0$ in 
$\overline{\Omega}\times[0,T_0],\ c(x,t)=1$ on 
$\overline{\Gamma_0}\times[0,T_0]$.  We say that 
potentials $A(x,t)=(A_0(x,t),A_1(x,t),...,A_n(x,t))$ and
$A'(x,t)=(A_0'(x,t),A_1'(x,t),...,A_n'(x,t))$
are gauge equivalent if there exists $c(x,t)\in G_0(\Omega\times[0,T_0])$ 
such that
\begin{eqnarray}                              \label{eq:1.6}
A_0'(x,t)=A_0(x,t)-ic^{-1}(x,t)\frac{\partial c}{\partial t},
\\
A_j'(x,t)=A_j(x,t)-ic^{-1}(x,t)\frac{\partial c}{\partial x_j},
\ \ 1\leq j\leq n.
\nonumber
\end{eqnarray}
Note that if $Lu=0$ and $u'=c(x,t)u$ then $L'u'=0$ where $L'$ is 
an operator of the form (\ref{eq:1.1}) with $A_j(x,t),\ 0\leq j\leq n,$
replaced by $A_j'(x,t),\ 0\leq j\leq n$.  We shall write for brevity
\[
L'=c\circ L.
\]
It is easy to show that if $\Lambda'$ is the D-to-N operator for $L'$
then $\Lambda'=\Lambda$ on $\Gamma_0\times(0,T_0)$,  i.e. all potentials
$A(x,t)$ in the same gauge equivalence class correspond to the same 
D-to-N operator on $\Gamma_0\times(0,T_0)$.  Note that if we consider
real-valued potentials only then the gauge group $G_0$ should be reduced 
to $c(x,t)$ such that $|c(x,t)|=1$.  If $\Omega$  is simply-connected then
any $c(x,t)\in G_0$ has a form $c(x,t)=e^{i\varphi(x,t)}$  where
$\varphi(x,t)\in C^\infty(\overline{\Omega}\times[0,T])$.  Also  if
coefficients of $L$ are independent of $t$ it is natural
 that the group $G_0$ consists of 
$c(x)$ independent of $t$.  Then $A_0'(x)=A_0(x)$ (see (\ref{eq:1.6}) ).
Denote
\[
T_*=\max_{x\in\overline{\Omega}} d(x,\Gamma_0),
\]
where $d(x,\Gamma_0)$ is the distance in $\overline{\Omega}$  with
respect to the metric $\|g^{jk}(x)\|^{-1}$ from $x\in \overline{\Omega}$
to $\Gamma_0$.
We shall assume $L$ and $\Gamma_0$ satisfy the 
 BLR-condition  (see [BLR92])
 for $t=T_{**}$.  This means roughly speaking that any
null-bicharacteristic of $L$ in $(\overline{\Omega}\times[0,T_{**}])\times
(\R^{n+1}\setminus\{0\})$ intersects $(\Gamma_0\times [0,T_{**}])\times
(\R^{n+1}\setminus\{0\})$.
It was proven in [BLR92]  that the BLR-condition implies that the 
(bounded)  map of $f\in H_0^1(\Gamma_0\times(0,T_0))$  to 
$(u(x,T_{**}),u_t(x,T_{**}))\in
H^1(\Omega)\times L^2(\Omega)$ is onto. Here $H_0^1(\Gamma_0\times (0,T_{**})$
is the subspace of $H^1(\partial\Omega\times(0,T_{**}))$  such that
$f|_{t=0}=0$  and $\mbox{supp\ }f\subset\overline{\Gamma_0}\times(0,T_{**}],
\ u(x,t)$ is the solution of (\ref{eq:1.1}), (\ref{eq:1.2}), (\ref{eq:1.3}).

The following theorem was proven in [E06]:

\begin{theorem}                              \label{theo:1.1}
Let $L$ and $L_0$ be two operators of the form (\ref{eq:1.1}) in domains
$\Omega$ and $\Omega_0$, respectively,
with coefficients $A(x,t), V(x,t)$  and $A_0(x,t), V_0(x,t)$ analytic 
in $t$ and real-valued.   Suppose  
$\Gamma_0\subset \partial\Omega\cap\partial\Omega_0$
and suppose that $L$ and $\Gamma_0$ satisfy the BLR-condition 
when $t=T_{**}$.   Suppose
that D-to-N operators $\Lambda$  and  $\Lambda_0$,  corresponding to $L$
and $L_0$, respectively,  are equal on $\Gamma_0\times(0,T_0)$ for all
smooth $f$ with supports on $\Gamma_0\times(0,T_0]$.  Let  $T_0>2T_*+T_{**}$.
Then there exists a diffeomorphism $\varphi$  of $\overline{\Omega}$ onto
$\overline{\Omega_0},\ \varphi=I$
on $\Gamma_0$,  and there exists a gauge transformation 
$c_0(x,t)\in G_0(\overline{\Omega}\times [0,T_0])$  such that
\[
c_0\circ \varphi^{-1}\circ L_0=L
\]
on $\Omega\times(0,T_0)$.
\end{theorem}

Denote by
$L^*$ the formally adjoint operator  to $L$.  
  Note that $L^*$ has the form (\ref{eq:1.1})
with $A_j(x,t),\ 0\leq j\leq n,
 V(x,t)$ replaced by $\overline{A_j(x,t)},\ 0\leq j\leq n,
 \overline{V(x,t)}$.  To prove Theorem \ref{theo:1.1} we need to know also 
the D-to-N operator $\Lambda_*$  corresponding to $L^*$.  If $L^*=L$  then
obviously $\Lambda_*=\Lambda$.
In the case when $A_0=0$  and potentials $A_j(x),\ 1\leq j\leq n,\ V(x)$ 
are independent of $t$  one can show that $\Lambda$  determines $\Lambda_*$
on $\Gamma_0\times(0,T_0)$  (c.f. [KL00] and \S 2 below)
  even when $A_j(x),\ 1\leq j\leq n,\ 
V(x)$ are complex-valued.
Therefore Theorem \ref{theo:1.1} holds in this case and gives a new proof
of the correponding result in [KL00], [KL2 97].
When $A_0=0,\ A_j(x),\ 1\leq j \leq n,\ V(x)$ are real-valued and 
independent of $t$,  i.e. in the self-adjoint case, 
the BLR-condition is not needed.  In this case Theorem \ref{theo:1.1}
is true with $T_{**}=0$.  This result was first obtained 
by BC (Boundary Control)  method (see [B97]) (see also  [B1 02], [B2 93], [KKL01].
[K93], [KK98]).  In [E1 06] we gave a new proof for
 time-independent self-adjoint case.  The proof in [E06]
is based on the new approach in [E1 06].  The inverse problems for
the wave equations with
 time-dependent
potentials in the case when $\Gamma_0=\partial\Omega$ was considered in
[St89], [RS91] (see also [I98]).

A crucial step of the proof of Theorem \ref{theo:1.1} uses the unique
continuation theorem by Tataru [T95].  This theorem requires that
$A_j(x,t),\ 0\leq j\leq n,
 V(x,t)$  depend analytically on $t$.

The proof of Theorem \ref{theo:1.1} consists of two steps : the local step
and the global step.  In the local step we recover the coefficients of $L$
(up to a diffeomorphism and a gauge transformation) in the domain 
$\Gamma_\delta\times[0,T_0]$  where $\Gamma$ is an open connected
subset of $\Gamma_0$ and $\Gamma_\delta$ in a small neighborhood of
$\Gamma$ in $\overline{\Omega}.$ 
The main novelty of the proof here is the study of the restrictions 
of the solutions of $Lu=0$ to the characteristic surfaces instead of
the restrictions to the hyperplanes $t=\mbox{constant}$  as in BC-method.

The main part of the global step is the following lemma that reduced the
inverse problem in the domain to the inverse problem in a smaller
domain (c.f. [KKL1 04]):
\begin{lemma}                                \label{lma:1.1}
Let $L^{(p)},p=1,2$ be two operators of the form (\ref{eq:1.1}) in
domains $\Omega_p,p=1,2,$ respectively,  satisfying the initial-boundary
conditions (\ref{eq:1.2}),  (\ref{eq:1.3}).  We assume that 
$\Gamma_0\subset\partial\Omega_1\cap\partial\Omega_2,\ 
\mbox{supp\ }f\subset \Gamma_0\times(0,T_0]$ and $\Lambda_1=\Lambda_2$
on $\Gamma_0\times (0,T_0)$ where $\Lambda_p$ are the D-to-N operators 
corresponding to $L^{(p)},\ p=1,2.$  Let $B\subset\Omega_1\cap\Omega_2$ be
such that the domains $\Omega_p\setminus\overline{B}$ are smooth,
$L^{(1)}=L^{(2)}$ in $\overline{B}$
 and
$S_1\stackrel{def}{=}\partial B\cap \partial\Omega_p\subset \Gamma_0,\ p=1,2.$
Let $\delta=\max_{x\in \overline{B}}d(x,\Gamma_0)$
where $d(x,\Gamma_0)$ is the distance in $\overline{B}$ 
from $x\in \overline{B}$
to $\Gamma_0$.  Denote by $\hat{\Lambda}_p$  the D-to-N operators 
corresponding to
$L^{(p)}$ in domains $(\Omega_p\setminus \overline{B})\times(\delta,T_0-\delta),
\ p=1,2$.  Let  $S_2=\partial B\setminus \overline{S_1}$ and let
$\Gamma_1=(\Gamma_0\setminus S_1)\cup S_2$.  Then $\hat{\Lambda}_1=
\hat{\Lambda}_2$ on $\Gamma_1\times(\delta,T_0-\delta)$.
\end{lemma}

\section{Hyperbolic systems with Yang-Mills  potentials and domains with obstacles.}
\label{section 2}
\init

Consider in $\Omega\times(0,T_0)$ a system of the form
(c.f. [E2 05])
\begin{eqnarray}                               \label{eq:2.1}
Lu\stackrel{def}{=}
\left(-i\frac{\partial u}{\partial t}I_m+A_0(x,t)\right)^2u(x,t)
\ \ \ \ \ \ \ \ \ \ \ \ \ \ \ \ \ \ \ \ \ \ \ \ \ \ \ \ \ \ 
\\
-\sum_{j,k=1}^n\frac{1}{\sqrt{g(x)}}
\left(-i\frac{\partial}{\partial x_j}I_m+A_j(x,t)
\right)
 \sqrt{g(x)}g^{jk}(x)
\left(-i\frac{\partial}{\partial x_k}I_m+A_k(x,t)\right)u
\nonumber
\\
-V(x,t)u=0,
\nonumber
\end{eqnarray}
where $u(x,t),\ A_j(x,t),\ 0\leq j\leq n,\ V(x,t)$ are $m\times m$ 
matrices,  $I_m$  is the identity $m\times m$ matrix.
Assume that the initial-boundary conditions   (\ref{eq:1.2}), (\ref{eq:1.3})
are satisfied.  Let $\Gamma_0\subset \partial\Omega$ and let
$G_0(\overline{\Omega}\times[0,T])$
be the gauge group of nonsingular $C^\infty\ \ m\times m$ 
matrices $C(x,t)$ in $\overline{\Omega}\times [0,T_0]$ analytic  in $t$
and such that $C(x,t)=I_m$ on $\Gamma_0\times [0,T_0]$.
Matrices $A(x,t)=(A_0(x,t),...,A_n(x,t)),\ V(x,t)$ are called
Yang-Mills potentials.  We say that $(A(x,t),V(x,t))$ and \\
$(A'(x,t),V'(x,t))$ are gauge equivalent if there exists 
$C(x,t)\in G_0(\Omega_0\times [0,T_0])$  such that 
\begin{eqnarray}                                 \label{eq:2.2}
A_0'(x,t)=C^{-1}(x,t)A_0(x,t)C(x,t)
-iC^{-1}(x,t)\frac{\partial C(x,t)}{\partial t},                                        
\\
A_j'(x,t)=C^{-1}A_j(x,t)C-iC^{-1}\frac{\partial C}{\partial x_j},\ \ 
1\leq j\leq n,
\nonumber
\\
V'(x,t)=C^{-1}V(x,t)C.
\nonumber
\end{eqnarray}

When we consider self-adjoint operators of the form (\ref{eq:2.1})
,
i.e.  when matrices $A_j(x,t),\ 0\leq j\leq n, \ V(x,t)$ are self-adjoint,
the group $G_0(\Omega\times [0,T_0])$ consists of unitary matrices $C(x,t)$.

A generalization of the proof of Theorem \ref{theo:1.1} 
leads to the following result   (c.f. [E2 05]:
\begin{theorem}                                \label{theo:2.1}
Theorem \ref{theo:1.1} holds for the equations of the form (\ref{eq:2.1})
with Yang-Mills potentials.
\end{theorem}
Consider now the system of the form (\ref{eq:2.1})  when the Yang-Mills 
potentials are independent of $t$ but not necessary self-adjoint matrices:
\begin{eqnarray}                               \label{eq:2.3}
Lu\stackrel{def}{=}
\left(-i\frac{\partial u}{\partial t}I_m+A_0(x)\right)^2u(x,t)
\ \ \ \ \ \ \ \ \ \ \ \ \ \ \ \ \ \ \ \ \ \ \ \ \ \ \ \ \ \ 
\\
-\sum_{j,k=1}^n\frac{1}{\sqrt{g(x)}}
\left(-i\frac{\partial}{\partial x_j}I_m+A_j(x)
\right)
 \sqrt{g(x)}g^{jk}(x)
\left(-i\frac{\partial}{\partial x_k}I_m+A_k(x)\right)u(x,t)
\nonumber
\\
-V(x)u(x,t)=0.
\nonumber
\end{eqnarray}
We also assume that $T_0=+\infty$,  i.e. (\ref{eq:2.3}) and
the boundary condition (\ref{eq:1.3}) hold for $t\in (0,+\infty)$.
Let $L^*$ be formally adjoint to $L$,  i.e.  when  
 $A_j(x),\ 0\leq j\leq n,\ V(x)$ are replaced
by the adjoint matrices $A_j^*(x),\ 0\leq j\leq n,\ V^*(x)$.

Consider the initial-boundary value problem adjoint to (\ref{eq:2.3}),
(\ref{eq:1.2}),   (\ref{eq:1.3})  on some interval $(0,T)$:
\begin{equation}                             \label{eq:2.4}
L^*v=0 \ \ \ \ \mbox{on\ \ \ } \Omega\times (0,T),
\end{equation}
\begin{equation}                             \label{eq:2.5}
v|_{t=T}=\frac{\partial v}{\partial t}|_{t=T}=0,
\ \ \ \ v|_{\partial\Omega\times (0,T)}=g,
\end{equation}
where $\mbox{supp\ }g\subset \Gamma_0\times(0,T]$.  Let $\Lambda^*$ 
be the D-to-N operator corresponding to (\ref{eq:2.4}),   (\ref{eq:2.5}).
We have
\[
0=(Lu,v)-(u,L^*v)=(\Lambda f,g)-(f,\Lambda^*g)
\]
for any smooth
$f$ and $g,\ \mbox{supp\ }f\subset\Gamma_0\times(0,T],\ 
\mbox{supp\ }g\subset \Gamma_0\times[0,T_0)$.  Therefore $\Lambda^*$
is an adjoint operator to $\Lambda$ and we can determine $\Lambda^*$ on 
$\Gamma_0\times[0,T)$ if we know $\Lambda$ on $\Gamma_0\times(0,T)$.  
Change in (\ref{eq:2.4}), (\ref{eq:2.5}) $t$ to $T-t$.  Then
we get an initial-boundary value problem
\begin{equation}                             \label{eq:2.6}
L_1^*w=0 \ \ \ \ \mbox{on\ \ \ } \Omega\times (0,T),
\end{equation}
\begin{equation}                             \label{eq:2.7}
w(x,0)=w_t(x,0)=0,
\ \ \ \ w|_{\partial\Omega\times (0,T)}=g_1(x,t),
\end{equation}
where $w(x,t)=v(x,T-t),\ g_1(x,t)=g(x,T-t),\ 0<t<T,\ L_1^*$ 
is obtained from $L^*$ by changing $A_0^*(x)$ to $-A_0^*(x)$.  It is
clear that the D-to-N operator $\Lambda_{1*}$ on 
$\Gamma_0\times(0,T)$ corresponding
to (\ref{eq:2.6}), (\ref{eq:2.7}) is determined by $\Lambda^*$.

Consider also the initial-boundary value problem
\begin{equation}                          \label{eq:2.8}
L^*u=0 \ \ \ \ \mbox{on\ \ \ } \Omega\times (0,T),
\end{equation}
\begin{equation}                             \label{eq:2.9}
u(x,0)=u_t(x,0)=0,
\ \ \ \ u|_{\partial\Omega\times (0,T)}=f(x,t).
\end{equation}
Denote by $\Lambda_*$ the D-to-N operator correspondin to
(\ref{eq:2.8}), (\ref{eq:2.9}). 
Here $T >0$ is arbitrary, i.e. 
(\ref{eq:2.6}), (\ref{eq:2.7}) and (\ref{eq:2.8}), (\ref{eq:2.9})
hold on $(0,+\infty)$.
We assume that $f(x,t)$ and $g_1(x,t)$ belong to 
  $C_0^\infty(\Gamma\times(0,+\infty))$.
Performing the Fourier-Laplace  transform in $t$ in
(\ref{eq:2.6}), (\ref{eq:2.7}) 
and in (\ref{eq:2.8}), (\ref{eq:2.9}) 
when $T=+\infty$  we get:
\begin{equation}                          \label{eq:2.10}
L^*(k)\tilde{u}(x,k)=0,     \ \ \ \ x\in\Omega,
\end{equation}
\begin{equation}                             \label{eq:2.11}
\tilde{u}(x,k)|_{\partial\Omega}=\tilde{f}(x,k),
\end{equation}
and
\begin{equation}                           \label{eq:2.12}
L^*(-k)\tilde{w}(x,k)=0,\ \ \ x\in\Omega,
\end{equation}
\begin{equation}                           \label{eq:2.13}
\tilde{w}(x,k)|_{\partial \Omega}=\tilde{g}_1(x,k),
\end{equation}
where $\tilde{u}(x,k),\tilde{w}(x,k)$ are analytic in $k$ for $\Im k <-C_0$
for some $C_0>0,\   L^*(k)$ is obtained from $L^*$ by replacing 
$-i\frac{\partial}{\partial t}$ by $k$.  Let $\Lambda_*(k)$ be
the D-to-N operator on $\Gamma$ corresponding
 to the boundary value problem (\ref{eq:2.10}), (\ref{eq:2.11}),
depending on parameter $k$.
Note that $\Lambda_*(k)$ is the Fourier-Laplace transform in $t$
of the D-to-N operator $\Lambda_{*}$ corresponding to
(\ref{eq:2.8}), (\ref{eq:2.9}) on $(0,+\infty)$.   Since $\Omega$  is 
a bounded domain $\Lambda_*(k)$ has an analytic continuation from  $\Re k\leq -C_0$
to $\C\setminus Z$  where $Z$ is a discrete set.  
Note
that the Fourier-Laplace transform of $\Lambda_{1*}$  is $\Lambda_*(-k)$.    
Since $\Lambda_*(k)$ is analytic in $\C\setminus Z$, 
$\Lambda_*(-k)$ determines $\Lambda_*(k)$.  Therefore when
$T_0=+\infty$ we get that the D-to-N  operator  $\Lambda$  on
$\Gamma_0\times(0,+\infty)$  determines the D-to-N operator $\Lambda_*$
on $\Gamma_0\times(0,+\infty)$.  Therefore the proof of Theorem \ref{theo:2.1}
applies  and we have the following result  (c.f. [KL1 00],  [KL2 97]:
\begin{theorem}                               \label{theo:2.2}
Let $L_p$ be two operators of the form (\ref{eq:2.3}) in domains 
$\Omega_p\times(0,+\infty),\ p=1,2$.   
Suppose $\Gamma_0\subset \partial\Omega_1\cap\partial\Omega_2$  and
$\Lambda_1=\Lambda_2$  on $\Gamma_0\times(0,+\infty)$  where $\Lambda_p$
are the D-to-N operators corresponding to $L_p,p=1,2$.
Suppose that $L_1$ and $\Gamma_0$ satisfy
the BLR-condition  
for some $t=T_{**}$.
Then there exists a diffeomorphism $y=\varphi(x)$ of  $\overline{\Omega_1}$
 onto $\overline{\Omega_2}$ 
and a gauge transformation $c_0(x)\in G_0(\overline{\Omega_1})$ such
that
\[
c_0\circ\varphi^{-1}\circ L_2=L_1, \ \ \ \ x\in\Omega.
\]
\end{theorem}
We do not assume here that $L_p,\ p=1,2,$  are formally self-adjoint.

Note that domains $\Omega$ can be multi-connected and 
$\Gamma_0\subset\partial\Omega$ can be not connected.
An important example of inverse problems with the boundary data prescribed 
on a part  of the boundary are the inverse problems in
domains with obstacles.  In this case
$\Omega=\Omega_0\setminus(\cup_{j=1}^r\overline{\Omega_j})$,  where
$\Omega_1,...,\Omega_r$  are nonintersecting domains inside $\Omega_0$,
called obstacles,
 $\Gamma_0=\partial\Omega_0$ and the zero Dirichlet boundary conditions
are prescribed on $\partial\Omega_j,\ 1\leq j\leq r$, i.e. we have
\begin{equation}                          \label{eq:2.14}
Lu=0 \ \ \ \mbox{on\ \ \ \ }\Omega\times(0,T_0),
\end{equation}  
\begin{eqnarray}                           \label{eq:2.15}
u(x,0)=u_t(x,0)=0\ \ \ \ \mbox{on}\ \ \ \Omega,
\\
u|_{\partial\Omega_0\times(0,T_0)}=f(x,t),\ \ \ \ 
u|_{\partial\Omega_j\times(0,T_0)}=0,\ \
1\leq j\leq r.
\nonumber
\end{eqnarray}
Unfortunately,  the BLR-condition is not satisfied for domains with more
than one smooth obstacle.  Therefore we shall assume that $L$ is a formally 
self-adjoint operator of the form (\ref{eq:2.3}),
i.e. when $A_j(x),\ 0\leq j\leq n,\ V(x)$ are self-adjoint matrices, and
initial-boundary conditions (\ref{eq:2.15}) are satisfied.  In this case
Theorem \ref{theo:2.1}  holds for any $T_0>2T_*$ and for any number of obstacles.

Finally  consider the following particular case: 
$T_0=+\infty,\ g^{jk}(x)=\delta_{jk},\ A_0(x)=0,\ A_j(x),\ 1\leq j\leq n,\ V(x)$
are self-adjoint. 
Making  the Fourier-Laplace transform 
in (\ref{eq:2.14})
 we get 
the Schr\"{o}dinger equation with Yang-Mills potentials in $\Omega$:
\begin{equation}                                    \label{eq:2.16}
\sum_{j=1}^n\left(-i\frac{\partial}{\partial x_j}I_m+A_j(x)\right)^2w(x)+
V(x)w(x)-k^2w(x)=0,
\end{equation}
where we omitted the dependence of $w$ on $k$ in (\ref{eq:2.16}).
When $m=1$ we have the Schr\"{o}dinger equation with electromagnetic potentials.
The boundary conditions for (\ref{eq:2.16}) 
have the form:
\begin{equation}                                \label{eq:2.17}
w|_{\partial\Omega_0}=h(x),\ \ \ w|_{\partial\Omega_j}=0,\ \ \ ,1\leq j\leq r.
\end{equation}
The D-to-N operator for (\ref{eq:2.16}), (\ref{eq:2.17}) has the form:
\begin{equation}                              \label{eq:2.18}
\Lambda(k)h=\left(\frac{\partial w}{\partial \nu}+i(A\cdot \nu)w\right)|_{\partial\Omega_0},
\end{equation}
where $\nu$ is the exterior unit normal vector to $\partial\Omega_0$.
Knowing the hyperbolic D-to-N operator for (\ref{eq:2.14}), (\ref{eq:2.15})
for the arbitrary $T_0>0$ we can find $\Lambda(k)$ for all $k\in \C\setminus Z$,
and vice versa.  Since $\Lambda(k)$ is analytic,  knowing $\Lambda(k)$ on
any interval $(k_0-\e,k_0+\e)$ of analyticity determines $\Lambda(k)$ for all
$k\in \C\setminus Z$.
Therefore  Theorem \ref{theo:2.1}  implies that $\Lambda(k)$ given on
$(k_0-\e,k_0+\e),\ k_0>0,\ \e>0$,  determines the location of all obstacles
$\Omega_j,\ 1\leq j\leq r,$  since the metric is fixed, and
determines potentials $A_j(x),\ 1\leq j\leq n,$
$V(x)$  in $\overline{\Omega}$ up to a gauge 
transformation $C(x)\in G_0(\overline{\Omega})$,  i.e  $C(x)=I_m$ on 
$\partial\Omega_0,\ C(x)$  is an unitary matrix in $\overline{\Omega}$.
 
The interest of considering multi-connected domains with obstacles was 
spurred by the Aharonov-Bohm effect.  It was shown by Aharonov and
Bohm [AB59] that the presence of distinct gauge equivalence classes of
potentials can be detected in an experiment and this phenomenon is
called the Aharonov-Bohm effect.  
As it was shown above the D-to-N $\Lambda(k)$
on $\partial\Omega_0$ given for all $k\in(k_0-\e,k_0+\e)$ allows to detect 
the gauge equivalent
class of Yang-Mills (or electroctromagnetic)  potentials.

\section{A geometric optic approach.}
\label{section 3}
\init

Consider the Schr\"{o}dinger equation with electromagnetic potentials
in the domain 
 $\Omega=\Omega_0\setminus (\cup_{j=1}^r\overline{\Omega_j})$
with obstacles, i.e.  consider  (\ref{eq:2.16}) when
$m=1$,  with boundary conditions (\ref{eq:2.17}).

Assume that the D-to-N operator $\Lambda(k)$ on $\partial\Omega_0$  (see
(\ref{eq:2.18}) )
is given for all $k\in \C\setminus K.$  Another approach
to the inverse problem for (\ref{eq:2.16}), (\ref{eq:2.17}) is based 
on geometric optics constructions and
the reduction to the integral geometry (tomography) problems.

We say that $\gamma=\gamma_1\cup\gamma_2\cup ... \cup\gamma_N$ is
a broken ray with legs $\gamma_1.\gamma_2,...,\gamma_N$ if 
$\gamma_k,\ 1\leq k\leq N$, are geodesics,
$\gamma$ starts at point $x_0\in\partial\Omega_0,\ \gamma$ has
$N-1$ nontangential points of reflection at the obstacles and
$\gamma$ ends at a point $x_N\in\partial\Omega_0$.
One can construct geometric optics solutions supported in a small neighborhood
of $\gamma$  (c.f. [E3 04], [E2 05])

Consider two Schr\"{o}dinger equations with electro-magnetic 
potentials \\ 
$A^{(p)}(x),\ V^{(p)}(x),\ p=1,2,$  with the Euclidian metric 
$g^{jk}=\delta_{jk}$
in a plane domain with convex obstacles.
Let $\Lambda_p(k)$ be the corresponding D-to-N operators,  $p=1,2$.

Using the geometric optics solutions one can prove that if the D-to-N 
operators 
are equal on $\partial\Omega_0$ then
\begin{equation}                           \label{eq:3.1}
 \exp(i\int_\gamma A^{(1)}(x)\cdot dx)
=\exp(i\int_\gamma A^{(2)}(x)\cdot dx),
\end{equation}
\begin{equation}                           \label{eq:3.2}
 \int_\gamma V^{(1)}(x)ds =\int_\gamma V^{(2)}(x)ds
\end{equation}
for any broken ray (c.f. [E3 04], [E2 05]).
The geometric optics construction and equalities 
(\ref{eq:3.1}),  (\ref{eq:3.2})
hold in any dimension  $n\geq 2$  and for any broken ray even when
the broken rays are passing through  generic caustics. 
Having (\ref{eq:3.1}),  (\ref{eq:3.2})  we  reduce the inverse problem for 
the Schr\"{o}dinger  equation to the inverse problem of the integral
geometry of broken rays,  i.e. the recovery of potentials from integrals
over broken rays.  This is a  difficult problem.

Some results in this direction were obtained in [E3 04] for $n=2$
under the geometric restriction that there is no trapped rays.  This
condition is not satisfied  when one has more than one smooth obstacle.
However,  there are piecewise smooth convex obstacles  that satisfy these
conditions.  
In this case it was shown in [E3 04]  that if (\ref{eq:3.1}),  (\ref{eq:3.2})
hold for all broken rays in $\Omega_0$ then $V^{(1)}=V^{(2)}$ and
$A^{(1)}$ and $A^{(2)}$ are gauge equivalent.
Despite that this approach is much more restrictive than
the hyperbolic equations  approach it has an advantage that it allows
to prove the stability results in some cases.
It also does not require the BLR-condition in the nonself-adjoint case.

Consider the following example:

Let $\Omega_1\subset\Omega_0$ be the only
 convex obstacle in $\Omega_0$ and let $f(x)$ be a smooth 
function in $\overline{\Omega_0}\setminus\Omega_1$.
It is well-known (c.f. [He80]) that if $\int_\gamma f(x)ds=0$ for all lines
$\gamma$ not intersecting $\overline{\Omega_1}$ then $f(x)=0$.
This problem is severly ill-posed.  If one uses the broken rays, i.e.
if one compute $\int_\gamma fds$ for all broken rays $\gamma$,
then the inverse problem is well-posed and there is a stability estimate.
More precisely, let $\gamma_{x,\theta}$ be the broken ray starting on
$\partial\Omega_0$  and ending at $x\in\overline{\Omega_0}\setminus\Omega_1$.
Here $\theta$ is the direction of the ray at the endpoint $x$.  We assume that
$w(x,\theta)$ is known when $x\in \partial\Omega_0,\ \forall\theta\in S^1$.
The following 
staility
estimate holds  (c.f. [E3 04] and [M77] in the case of no obstacles) :
\[
\int_{\Omega_0\setminus\Omega_1}|f(x)|^2dx\leq
 C\int_0^{l_0}\int_{S^1}
\left(\left|\frac{\partial w(x(s),\theta)}{\partial s}\right|
+\left|\frac{\partial w}{\partial \theta}\right|^2\right)dsd\theta,
\]
where $x=x(s)$ is the equation of $\partial\Omega_0$,
$l_0$ is the arclength of $\partial\Omega_0$.

\end{document}